\documentclass{amsart}
\usepackage{cases}
\usepackage{latexsym}
\usepackage[arrow,matrix]{xy}
\usepackage{stmaryrd}
\usepackage{amsfonts}
\usepackage{amsmath,amssymb,amscd,bbm,amsthm,mathrsfs,dsfont}
\usepackage{fancyhdr}
\usepackage{amsxtra,ifthen}
\usepackage{verbatim}

\numberwithin{equation}{section}

\theoremstyle{plain}
\newtheorem{prop}{Proposition}[section]
\newtheorem{thm}[prop]{Theorem}
\newtheorem{cor}[prop]{Corollary}
\newtheorem{lem}[prop]{Lemma}
\newtheorem{dfn}[prop]{Definition}

\theoremstyle{definition}

\newtheorem{remark}[prop]{Remark}

\newcommand{\lam}{\lambda}

\DeclareMathOperator{\rad}{rad}
\DeclareMathOperator{\rank}{rank}

\begin{document}

\title{Radicals of symmetric cellular algebras}
\thanks {This work is partially supported by the Research Fund of
Doctor Program of Higher Education, Ministry of Education of
China.}
\maketitle

\centerline{Yanbo Li}
\begin{center}Department of Information and Computing Sciences,
\\Northeastern University at Qinhuangdao; \\Qinhuangdao, 066004, P.R. China \\
School of Mathematics Sciences, Beijing Normal University;\\
Beijing, 100875, P.R. China\\
\texttt{E-mail: liyanbo707@163.com}
\end{center}

\begin{abstract}
Using a slightly weaker definition of cellular algebra, due to Goodman (\cite{G2} Definition 2.9), we prove that for a symmetric cellular algebra, the dual basis of a cellular basis is again cellular. Then a nilpotent ideal is constructed for a symmetric cellular algebra. The ideal connects the radicals of cell modules with the radical of the algebra. It also reveals some information on the dimensions of simple modules. As a by-product, we obtain some equivalent conditions for a finite dimensional symmetric cellular algebra to be semisimple.
\end{abstract}

\medskip {\small {\it 2000 AMS Classification}: 16G30, 16N20

}

\bigskip

\section{\bf Introduction}
\label{xxsec1}

Cellular algebras were introduced by Graham and Lehrer \cite{GL} in
1996, motivated by previous work of Kazhdan and Lusztig \cite{KL}.
They were defined by a so-called cellular basis with some nice
properties. The theory of cellular algebras provides a systematic
framework for studying the representation theory of non-semisimple
algebras which are deformations of semisimple ones. One can
parameterize simple modules for a finite dimensional cellular
algebra by methods in linear algebra. Many classes of algebras from
mathematics and physics are found to be cellular, including Hecke
algebras of finite type, Ariki-Koike algebras, $q$-Schur algebras,
Brauer algebras, Temperley-Lieb algebras, cyclotomic Temperley-Lieb
algebras, Jones algebras, partition algebras, Birman-Wenzl algebras
and so on, we refer the reader to \cite{G, GL, RX, Xi1, Xi2} for details.

An equivalent basis-free definition of cellular algebras was given by Koenig and Xi
\cite{KX1}, which is useful in dealing with structural problems.
Using this definition, in
\cite{KX5}, Koenig and Xi made explicit an inductive construction of
cellular algebras called inflation, which produces all cellular
algebras. In \cite{KX7}, Brauer algebras were shown to be iterated
inflations of group algebras of symmetric groups and then more
information about these algebras was found.

There are some generalizations of cellular algebras, we refer the reader to \cite{DR, GRM, GRM2, WB} for details. Recently, Koenig and
Xi \cite{KX8} introduced affine cellular algebras which contain
cellular algebras as special cases. Affine Hecke algebras of type A
and infinite dimensional diagram algebras like the affine
Temperley-Lieb algebras are affine cellular.

\smallskip

It is an open problem to find explicit formulas for the dimensions
of simple modules of a cellular algebra. By the theory of cellular
algebras, this is equivalent to determine the dimensions of the
radicals of bilinear forms associated with cell modules. In \cite{LZ}, for a
quasi-hereditary cellular algebra, Lehrer and Zhang found that
the radicals of bilinear forms are related to the
radical of the algebra. This leads us to studying the
radical of a cellular algebra. However, we have no idea for dealing with general cellular algebras now. We will do some work on the radicals of
{\em symmetric} cellular algebras in this paper. Note that Hecke algebras of finite types, Ariki-Koike algebras
over any ring containing inverses of the parameters, Khovanov's diagram algebras are all symmetric cellular algebras.
The trivial extension of a cellular algebra is also a symmetric cellular algebra. For details, see \cite{BS}, \cite{MM}, \cite{XX}.

Throughout this paper, we will adopt a slightly weaker definition of cellular algebra due to Goodman (\cite{G2} Definition 2.9). It is helpful to note that the results of \cite{GL} remained valid with his weaker axiom. In case $2$ is invertible, these two definitions are equivalent.

We begin with recalling definitions and some well-known results of symmetric algebras and cellular algebras in Section 2. Then in Section 3, we prove that for a symmetric cellular algebra, the dual basis of a cellular basis is again cellular. In Section 4, a nilpotent ideal of a symmetric cellular algebra is constructed. This ideal connects the radicals of cell modules with the radical of the algebra and also reveals some information on the dimensions of simple modules. As a by-product, in Section 5, we obtain some equivalent conditions for a finite dimensional symmetric cellular algebra to be semisimple.

\bigskip

\section{\bf Preliminaries}
\label{xxsec2}
In this section, we start with the definitions of symmetric algebras and cellular algebras (a slightly weaker version due to Goodman) and then recall some well-known results about them.

Let $R$ be a commutative ring with identity and $A$ an associative
$R$-algebra. As an $R$-module, $A$ is finitely generated and free. Suppose that there exists an $R$-bilinear map
$f:A\times A\rightarrow R$. We say that $f$ is
 non-degenerate if the determinant of the matrix
$(f(a_{i},a_{j}))_{a_{i},a_{j}\in B}$ is a unit in $R$ for some
$R$-basis $B$ of $A$. We say $f$ is associative if $f(ab,c)=f(a,bc)$
for all $a,b,c\in A$, and symmetric if $f(a,b)=f(b,a)$ for all
$a,b\in A$.
\begin{dfn}\label{2.1}
An $R$-algebra $A$ is called symmetric if there is a non-degenerate
associative symmetric bilinear form $f$ on $A$. Define an $R$-linear map $\tau: A\rightarrow R$ by $\tau(a)=f(a,1)$.
We call $\tau$ a symmetrizing trace.
\end{dfn}

Let $A$ be a symmetric algebra with a basis $B=\{a_{i}\mid
i=1,\ldots,n\}$ and $\tau$ a symmetrizing trace. Denote by
$D=\{D_{i}\mid i=i,\ldots,n\}$ the basis determined by the
requirement that $\tau(D_{j}a_{i})=\delta_{ij}$ for all $i,
j=1,\ldots,n$. We will call $D$ the dual basis of $B$. For arbitrary $1\leq i,j \leq n$, write
$a_{i}a_{j}=\sum\limits_{k}r_{ijk}a_{k}$, where $r_{ijk}\in R$. Fixing
a symmetrizing trace $\tau$ for $A$, then we have the following lemma.
\begin{lem}\label{2.2}
Let $A$ be a symmetric $R$-algebra with a basis $B$ and the dual
basis $D$. Then the following hold:
$$a_{i}D_{j}=\sum_{k}r_{kij}D_{k};\,\,\,\,\,D_{i}a_{j}=\sum_{k}r_{jki}D_{k}.$$
\end{lem}

\begin{proof}
We only prove the first equation. The other one
is proved similarly.

Suppose that $a_{i}D_{j}=\sum\limits_{k}r_{k}D_{k}$, where $r_{k}\in R$ for $k=1,\cdots,n$. Left multiply by
$a_{k_{0}}$ on both sides of the equation and then apply $\tau$, we get
$\tau(a_{k_{0}}a_{i}D_{j})=r_{k_{0}}$. Clearly, $\tau(a_{k_{0}}a_{i}D_{j})=r_{k_{0},i,j}$. This implies that
$r_{k_{0}}=r_{k_{0},i,j}$.
\end{proof}

Given a symmetric algebra, it is natural to consider the relation between two dual bases determined by two different symmetrizing traces. For this we have the following lemma.

\begin{lem}\label{2.3}
Suppose that $A$ is a symmetric $R$-algebra with a
basis $B=\{a_{i}\mid i=1, \cdots, n\}$.  Let $\tau, \tau'$ be two
symmetrizing traces. Denote by $\{D_{i}\mid i=1, \cdots, n\}$ the
dual basis of $B$ determined by $\tau$ and $\{D_{i}'\mid i=1,
\cdots, n\}$ the dual basis determined by $\tau'$. Then for
$1\leq i \leq n$, we have
$$D_{i}'=\sum_{j=1}^{n}\tau(a_{j}D_{i}')D_{j}.$$
\end{lem}

\begin{proof}
It is proved by a similar method as in Lemma \ref{2.2}.
\end{proof}

\smallskip

Graham and Lehrer introduced the so-called cellular algebras in \cite{GL} , then Goodman weakened the definition in \cite{G2}. We will adopt Goodman's definition throughout this paper.

\begin{dfn} {\rm (\cite{G2})}\label{2.4}
Let $R$ be a commutative ring with identity. An
associative unital $R$-algebra is called a cellular algebra with
cell datum $(\Lambda, M, C, i)$ if the following conditions are
satisfied:

{\rm(C1)} The finite set $\Lambda$ is a poset. Associated with each
$\lam\in\Lambda$, there is a finite set $M(\lam)$. The algebra $A$
has an $R$-basis $\{C_{S,T}^\lam \mid S,T\in
M(\lam),\lam\in\Lambda\}$.

{\rm(C2)} The map $i$ is an $R$-linear anti-automorphism of $A$ with
$i^{2}=id$ and $$i(C_{S,T}^\lam)\equiv C_{T,S}^\lam  \,\,\,\,(\rm {mod}\,\,\,
 A(<\lam))$$ for all $\lam\in\Lambda$ and $S,T\in M(\lam)$, where $A(<\lam)$ is the $R$-submodule of $A$ generated by
$\{C_{S^{''},T^{''}}^\mu \mid S^{''},T^{''}\in M(\mu),\mu<\lam\}$.

{\rm(C3)} If $\lam\in\Lambda$ and $S,T\in M(\lam)$, then for any element
$a\in A$, we have\\
$$aC_{S,T}^\lam\equiv\sum_{S^{'}\in
M(\lam)}r_{a}(S',S)C_{S^{'},T}^{\lam} \,\,\,\,(\rm {mod}\,\,\,
 A(<\lam)),$$ where $r_{a}(S^{'},S)\in R$ is independent of $T$.

Apply $i$ to the equation in {\rm(C3)}, we obtain

{\rm(C3$'$)} $C_{T,S}^\lam i(a)\equiv\sum\limits_{S^{'}\in
M(\lam)}r_{a}(S^{'},S)C_{T,S^{'}}^{\lam} \,\,\,\,(\rm mod
\,\,\,A(<\lam)).$
\end{dfn}

\begin{remark}
Graham and Lehrer's original definition in \cite{GL} requires that $i(C_{S,T}^\lam)=C_{T,S}^\lam$ for all $\lam\in\Lambda$ and $S,T\in M(\lam)$. But Goodman pointed out that the results of \cite{GL} remained valid with his weaker axiom. In case $2\in R$ is invertible, these two definitions are equivalent.
\end{remark}

It is easy to check the following lemma by Definition \ref{2.4}.

\begin{lem} {\rm (\cite{GL})}
Let $\lam\in\Lambda$ and $a\in A$. Then for arbitrary elements
$S,T,U,V\in M(\lam)$, we have $$C_{S,T}^\lam aC_{U,V}^\lam \equiv
\Phi_{a}(T,U)C_{S,V}^\lam\,\,\,\, (\rm mod\,\,\, A(<\lam)),$$ where
$\Phi_{a}(T,U)\in R$ depends only on $a$, $T$ and $U$.
\end{lem}

We often omit the index $a$ when $a=1$, that is, writing
$\Phi_{1}(T,U)$ as $\Phi(T,U)$.

Let us recall the definition of cell modules now.
\begin{dfn} {\rm (\cite{GL})}\label{2.7}
Let $A$ be a cellular algebra with cell datum $(\Lambda, M, C, i)$. For each $\lam\in\Lambda$, define the left $A$-module
$W(\lam)$ as follows: $W(\lam)$ is a free $R$-module with basis
$\{C_{S}\mid S\in M(\lam)\}$ and $A$-action
defined by\\
$$aC_{S}=\sum_{S^{'}\in M(\lam)}r_{a}(S^{'},S)C_{S^{'}}
\,\,\,\,(a\in A,S\in M(\lam)),$$ where $r_{a}(S^{'},S)$ is the
element of $R$ defined in Definition {\rm\ref{2.4}} {\rm(C3)}.
\end{dfn}

Note that $W(\lam)$ may be thought of as a right $A$-module via \\
$$C_{S}a=\sum_{S^{'}\in M(\lam)}r_{i(a)}(S^{'},S)C_{S^{'}}
\,\,\,\,(a\in A,S\in M(\lam)).$$ We will denote this right $A$-module by $i(W(\lam))$.

\begin{lem}{\rm(\cite{GL})}\label{2.8}
There is a natural isomorphism of $R$-modules
$$C^{\lam}:W(\lam)\otimes_{R}i(W(\lam))\rightarrow R{\rm -span}\{C_{S,T}^{\lam}\mid S,T\in M(\lam)\},$$
defined by $(C_{S},C_{T})\rightarrow C_{S,T}^{\lam}$.
\end{lem}

For a cell module $W(\lam)$, define a bilinear form $\Phi
_{\lam}:\,\,W(\lam)\times W(\lam)\longrightarrow R$ by $\Phi
_{\lam}(C_{S},C_{T})=\Phi(S,T)$. It plays an important role for studying the structure of $W(\lam)$. It is easy to check that $\Phi(T,U)=\Phi(U,T)$ for arbitrary $T,U\in M(\lam)$.

Define
$$\rad\lam:= \{x\in W(\lam)\mid \Phi_{\lam}(x,y)=0
\,\,\,\text{for all} \,\,\,y\in W(\lam)\}.$$ If $\Phi
_{\lam}\neq 0$, then $\rad\lam$ is the radical of the $A$-module
$W(\lam)$. Moreover, if $\lam$ is a maximal element in $\Lambda$, then $\rad\lam=0$.

The following results were proved by Graham and Lehrer in \cite{GL}.

\begin{thm} {\rm\cite{GL}}
Let $K$ be a field and $A$ a finite dimensional cellular algebra.  For any
$\lam\in\Lambda$, denote the $A$-module $W(\lam)/\rad \lam$ by $L_{\lam}$. Let
$\Lambda_{0}=\{\lam\in\Lambda\mid \Phi_{\lam}\neq 0\}$. Then
$\{L_{\lam}\mid \lam\in\Lambda_{0}\}$ is a complete set of
{\rm (}representative of equivalence classes of {\rm )} absolutely simple
$A$-modules.
\end{thm}

\begin{thm}{\rm(\cite{GL})}
\label{glthm}Let $K$ be a field and $A$ a cellular $K$-algebra. Then the
following are equivalent.\\
{\rm(1)} The algebra $A$ is semisimple.\\
{\rm(2)} The nonzero cell representations $W(\lam)$ are irreducible and
pairwise inequivalent.\\
{\rm(3)} The form $\Phi_{\lam}$ is non-degenerate (i.e. $\rad\lam=0$)
for each $\lam\in\Lambda$.
\end{thm}

For any $\lam\in\Lambda$, fix an order on $M(\lam)$ and let
$M(\lam)=\{S_{1},S_{2},\cdots,S_{n_{\lam}}\}$, where $n_{\lam}$ is
the number of elements in $M(\lam)$,
the matrix $G(\lam)=(\Phi(S_{i},S_{j}))_{1\leq i,j\leq n_{\lam}}$ is called  Gram matrix.
It is easy to know that all the determinants of $G(\lam)$ defined with different order on $M(\lam)$ are the same.
By the definition of $G(\lam)$ and $\rad\lam$, for a finite dimensional
cellular algebra $A$, it is clear that if $\Phi_{\lam}\neq 0$, then $\dim_{K}L_{\lam}=\rank G(\lam)$.

\bigskip

\section{\bf Symmetric cellular algebras}

In this section, we prove that for a symmetric cellular algebra, the dual basis of a cellular basis is again cellular.

Let $A$ be a symmetric cellular algebra with a cell datum $(\Lambda,
M, C, i)$. Denote the dual basis by $D=\{D_{S,T}^\lam \mid
S,T\in M(\lam),\lam\in\Lambda\}$ throughout, which satisfies
$$
\tau(C_{S,T}^{\lam}D_{U,V}^{\mu})=\delta_{\lam\mu}\delta_{SV}\delta_{TU}.
$$
For any $\lam, \mu\in \Lambda$, $S,T\in M(\lam)$, $U,V\in M(\mu)$,
write
$$C_{S,T}^{\lam}C_{U,V}^{\mu}=\sum\limits_{\epsilon\in\Lambda,X,Y\in M(\epsilon)}
r_{(S,T,\lam),(U,V,\mu),(X,Y,\epsilon)}C_{X,Y}^{\epsilon}.$$
A lemma which we now prove plays an important role throughout this paper.
\begin{lem}\label{2.14}
Let $A$ be a symmetric cellular algebra with a cell datum $(\Lambda,
M, C, i)$ and $\tau$ a given symmetrizing trace. For arbitrary $\lam,\mu\in\Lambda$ and $S,T,P,Q\in M(\lam)$, $U,V\in M(\mu)$, the following hold:\\
{\rm(1)}\,\,\,\,$D_{U,V}^{\mu}C_{S,T}^{\lam}=\sum\limits_{\epsilon\in
\Lambda, X,Y\in
M(\epsilon)}r_{(S,T,\lam),(Y,X,\epsilon),(V,U,\mu)}D_{X,Y}^{\epsilon}.$\\
{\rm(2)}\,\,\,\,$C_{S,T}^{\lam}D_{U,V}^{\mu}=\sum\limits_{\epsilon\in
\Lambda, X,Y\in
M(\epsilon)}r_{(Y,X,\epsilon),(S,T,\lam),(V,U,\mu)}D_{X,Y}^{\epsilon}.$\\
{\rm(3)}\,\,\,\,$C_{S,T}^{\lam}D_{T,Q}^{\lam}=C_{S,P}^{\lam}D_{P,Q}^{\lam}.$\\
{\rm(4)}\,\,\,\,$D_{T,S}^{\lam}C_{S,Q}^{\lam}=D_{T,P}^{\lam}C_{P,Q}^{\lam}.$\\
{\rm(5)}\,\,\,\,$C_{S,T}^{\lam}D_{P,Q}^{\lam}=0\,\, if \,\,T\neq P.$\\
{\rm(6)}\,\,\,\,$D_{P,Q}^{\lam}C_{S,T}^{\lam}=0\,\, if \,\,Q\neq S.$\\
{\rm(7)}\,\,\,\,$C_{S,T}^{\lam}D_{U,V}^{\mu}=0 \,\,\,\,if\,\,\, \mu\nleq \lam.$\\
{\rm(8)}\,\,\,\,$D_{U,V}^{\mu}C_{S,T}^{\lam}=0 \,\,\,\,if \,\,\,\mu\nleq
\lam.$
\end{lem}
\begin{proof} (1), (2) are corollaries of Lemma \ref {2.2}. The equations
(5), (6), (7), (8) are corollaries of (1) and (2). We now prove (3).

By (2), we have
$$C_{S,T}^\lam D_{T,Q}^\lam=\sum_{\epsilon\in
\Lambda, X,Y\in
M(\epsilon)}r_{(Y,X,\epsilon),(S,T,\lam),(Q,T,\lam)}D_{X,Y}^{\epsilon}$$
$$C_{S,P}^\lam
D_{P,S}^\lam=\sum_{\epsilon\in \Lambda, X,Y\in
M(\epsilon)}r_{(Y,X,\epsilon),(S,P,\lam),(Q,P,\lam)}D_{X,Y}^{\epsilon}.$$
On the other hand, by (C3) of Definition \ref{2.4} we also have
$$r_{(Y,X,\epsilon),(S,T,\lam),(Q,T,\lam)}=r_{(Y,X,\epsilon),(S,P,\lam),(Q,P,\lam)}$$
for all $\epsilon\in \Lambda$ and $X,Y\in M(\epsilon)$. This
completes the proof of (3).

(4) is proved similarly.
\end{proof}

\begin{lem}\label{2.15}
Let $A$ be a symmetric cellular algebra with a
cell datum $(\Lambda, M, C, i)$. Then the dual basis $D=\{D_{S,T}^\lam \mid
S,T\in M(\lam),\lam\in\Lambda\}$ is again a cellular basis of $A$ with respect to the opposite order on $\Lambda$.
\end{lem}

\begin{proof} Clearly, we only need to consider (C2) and (C3) of Definition \ref{2.4}. Now we proceed in two steps.\\

{\it Step 1.} (C2) holds.

Let $i(D_{S,T}^{\lam})=\sum\limits_{\epsilon\in\Lambda, X,Y\in M(\epsilon)}r_{X,Y,\epsilon}D_{X,Y}^{\epsilon}$ with $r_{X,Y,\epsilon}\in R$. If there exists $\eta\ngeq\lam$ such that $r_{P,Q,\eta}\neq 0$ for some $P,Q\in M(\eta)$. Then $\tau(i(D_{S,T}^{\lam})C_{Q,P}^\eta)=r_{P,Q,\eta}\neq 0$. This implies that $i(D_{S,T}^{\lam})C_{Q,P}^\eta\neq 0$. Thus $C_{P,Q}^{\eta}D_{S,T}^{\lam}\neq 0$. But we know $\eta\ngeq\lam$, then by Lemma \ref{2.14} (7), $C_{P,Q}^{\eta}D_{S,T}^{\lam}=0$, a contradiction. This implies that $$i(D_{S,T}^{\lam})\equiv \sum\limits_{X,Y\in M(\lam)}r_{X,Y,\lam}D_{X,Y}^{\lam}\,\,\,(\mod A_{D}(>\lam)).$$
Now assume $r_{U,V,\lam}\neq 0$. Then $i(D_{S,T}^{\lam})C_{V,U}^{\lam}\neq 0$, hence $C_{U,V}^{\lam}D_{S,T}^{\lam}\neq 0$. By Lemma \ref{2.14} (5), $V=S$. We can get $U=T$ similarly.\\

{\it Step 2.} (C3) holds.

For arbitrary $C_{S,T}^{\lam}$, by Lemma \ref{2.14} (2), we have
$$C_{S,T}^{\lam}D_{U,V}^{\mu}=\sum_{\epsilon\in\Lambda, X,Y\in
M(\epsilon)}r_{(Y,X,\epsilon),(S,T,\lam),(V,U,\mu)}D_{X,Y}^{\epsilon}.$$
By (C3) of Definition \ref{2.4}, if $\epsilon <\mu$, then $r_{(Y,X,\epsilon),(S,T,\lam),(V,U,\mu)}=0$. Therefore,
$$C_{S,T}^{\lam}D_{U,V}^{\mu}\equiv\sum_{X,Y\in
M(\mu)}r_{(Y,X,\mu),(S,T,\lam),(V,U,\mu)}D_{X,Y}^{\mu}\,\,\,\,\,(\mod A_{D}(>\mu)),$$
where $A_{D}(>\mu)$ is the $R$-submodule of $A$ generated by
$$\{D_{S^{''},T^{''}}^\eta \mid S^{''},T^{''}\in M(\lam),\eta>\mu\}.$$
By (C3$'$) of Definition \ref{2.4}, if $Y\neq V$, then $r_{(Y,X,\mu),(S,T,\lam),(V,U,\mu)}=0$. So
$$C_{S,T}^{\lam}D_{U,V}^{\mu}\equiv\sum_{X\in M(\mu)}r_{(V,X,\mu),(S,T,\lam),(V,U,\mu)}D_{X,V}^{\mu}\,\,\,(\mod A_{D}(>\mu)).$$
Clearly, for
arbitrary $X\in M(\mu)$, we have
$$r_{(V,X,\mu),(S,T,\lam),(V,U,\mu)}=r_{C_{T,S}^{\lam}}(U,X)$$ and which is
independent of $V$. Since $C_{S,T}^{\lam}$ is arbitrary, then
$$aD_{U,V}^{\mu}\equiv \sum_{U'\in
M(\mu)}r_{i(a)}(U,U')D_{U',V}^{\mu}\,\,\,\,\,(\mod
A_{D}(>\mu))$$ for any $a\in A$. By Definition \ref{2.4}, $r_{i(a)}(U,U')$ is independent of $V$.
\end{proof}

\begin{remark}
Using the original definition of cellular algebras, Graham proved in \cite{G3} the dual basis of a cellular basis is again cellular in the case when $\tau(a)=\tau(i(a))$, for all $a\in A$.
\end{remark}

Since the dual basis is again cellular,
for arbitrary elements
$S,T,U,V\in M(\lam)$, it is clear that $$D_{S,T}^\lam D_{U,V}^\lam \equiv
\Psi(T,U)D_{S,V}^\lam\,\,\,\, (\rm mod\,\,\, A(>\lam)),$$ where
$\Psi(T,U)\in R$ depends only on $T$ and $U$. Then we also have Gram matrices $G'(\lam)$ defined by the dual basis. Now it is natural to consider the problem what is the relation between $G(\lam)$ and $G'(\lam)$. To study this, we need the following lemma.

\begin{lem}
\label{2.17}Let $A$ be a symmetric cellular algebra with cell datum $(\Lambda, M, C, i)$. For every
$\lam\in\Lambda$ and $S,T,U,V,P\in M(\lam)$, we have
$$C_{S,T}^{\lam}D_{T,U}^{\lam}C_{U,V}^{\lam}D_{V,P}^{\lam}=\sum_{Y\in
M(\lam)}\Phi(Y,V)\Psi(Y,V)C_{S,T}^{\lam}D_{T,P}^{\lam}.$$
\end{lem}
\begin{proof}
By Lemma \ref{2.14} (1), we have
\begin{eqnarray*}
& &
C_{S,T}^{\lam}D_{T,U}^{\lam}C_{U,V}^{\lam}D_{V,P}^{\lam}
=C_{S,T}^{\lam}(D_{T,U}^{\lam}C_{U,V}^{\lam})D_{V,P}^{\lam}\\&=&\sum_{\epsilon\in\Lambda,X,Y\in
M(\epsilon)}r_{(U,V,\lam),(Y,X,\epsilon),(U,T,\lam)}C_{S,T}^{\lam}D_{X,Y}^{\epsilon}D_{V,P}^{\lam}.
\end{eqnarray*}
If $\varepsilon>\lam$, then by Lemma \ref{2.14} (7), $C_{S,T}^{\lam}D_{X,Y}^{\epsilon}=0$; if $\varepsilon<\lam$, by Definition \ref{2.4} (C3), $r_{(U,V,\lam),(Y,X,\epsilon),(U,T,\lam)}=0$. This implies that
\begin{eqnarray*}
& &
\sum_{\epsilon\in\Lambda,X,Y\in
M(\epsilon)}r_{(U,V,\lam),(Y,X,\epsilon),(U,T,\lam)}C_{S,T}^{\lam}D_{X,Y}^{\epsilon}D_{V,P}^{\lam}\\&=&\sum_{X,Y\in
M(\lam)}r_{(U,V,\lam),(Y,X,\lam),(U,T,\lam)}C_{S,T}^{\lam}D_{X,Y}^{\lam}D_{V,P}^{\lam}.
\end{eqnarray*}
By Definition \ref{2.4} (C3), if $X\neq T$, then $r_{(U,V,\lam),(Y,X,\lam),(U,T,\lam)}=0$. Hence,
\begin{eqnarray*}
& &
\sum_{X,Y\in
M(\lam)}r_{(U,V,\lam),(Y,X,\lam),(U,T,\lam)}C_{S,T}^{\lam}D_{X,Y}^{\lam}D_{V,P}^{\lam}\\&=&\sum_{Y\in
M(\lam)}r_{(U,V,\lam),(Y,T,\lam),(U,T,\lam)}C_{S,T}^{\lam}D_{T,Y}^{\lam}D_{V,P}^{\lam}.
\end{eqnarray*}
Note that $$D_{T,Y}^{\lam}D_{V,P}^{\lam}\equiv\Psi(Y,V)D_{T,P}^{\lam} \,\,\,\,\,\,(\mod A_{D}(>\lam)).$$
Moreover, by Lemma \ref{2.14} (7), if $\epsilon>\lam$, then $C_{S,T}^{\lam}D_{X,Y}^{\epsilon}=0$. Thus
$$\sum_{Y\in
M(\lam)}r_{(U,V,\lam),(Y,T,\lam),(U,T,\lam)}C_{S,T}^{\lam}D_{T,Y}^{\lam}D_{V,P}^{\lam}=\sum\limits_{Y\in
M(\lam)}\Phi(Y,V)\Psi(Y,V)C_{S,T}^{\lam}D_{T,P}^{\lam}.$$
This completes the proof.
\end{proof}

By Lemma \ref{2.14}, $C_{U,V}^{\lam}D_{V,P}^{\lam}$ is
independent of $V$, so is $\sum\limits_{Y\in
M(\lam)}\Phi(Y,V)\Psi(Y,V)$. Then for any $\lam\in\Lambda$, we can
define a constant $k_{\lam,\tau}$ as follows.
\begin{dfn}\label{2.18} Keep the notation above.
For $\lam\in\Lambda$, take an arbitrary $V\in M(\lam)$. Define
$$k_{\lam, \tau}=\sum\limits_{X\in M(\lam)}\Phi(X,V)\Psi(X,V).$$
\end{dfn}

Note that $\{k_{\lam, \tau}\mid\lam\in\Lambda\}$ is not independent
of the choice of symmetrizing trace. Fixing a symmetrizing trace $\tau$,
we often write $k_{\lam, \tau}$ as $k_{\lam}$. The following lemma reveals the relation among $G(\lam)$, $G'(\lam)$ and $k_{\lam}$.
\begin{lem}\label{2.19}
Let $A$ be a symmetric cellular algebra with cell datum $(\Lambda, M, C, i)$. For any $\lam\in\Lambda$, fix an order on the set
$M(\lam)$. Then $G(\lam)G'(\lam)=k_{\lam}E$, where
$E$ is the identity matrix.
\end{lem}
\begin{proof}
For an arbitrary $\lam\in\Lambda$, according to the definition of $G(\lam)$, $G'(\lam)$ and $k_{\lam}$, we only
need to show that $\sum\limits_{Y\in M(\lam)}\Phi(Y,U)\Psi(Y,V)=0$ for
arbitrary $U,V\in M(\lam)$ with $U\neq V$.

In fact, on one hand, for arbitrary $S\in M(\lam)$, by Lemma
\ref{2.14} (5), $U\neq V$ implies that $C_{S,U}^{\lam}D_{V,S}^{\lam}=0$. Then
$C_{S,U}^{\lam}D_{U,S}^{\lam}C_{S,U}^{\lam}D_{V,S}^{\lam}=0$.

On the other hand, by a similar method as in the proof of Lemma \ref{2.17},
\begin{eqnarray*}
C_{S,U}^{\lam}D_{U,S}^{\lam}C_{S,U}^{\lam}D_{V,S}^{\lam}&=&\sum_{\epsilon\in\Lambda,X,Y\in
M(\epsilon)}r_{(S,U,\lam),(Y,X,\epsilon),(S,U,\lam)}C_{S,U}^{\lam}D_{X,Y}^{\epsilon}D_{V,S}^{\lam}\\
&=&\sum_{Y\in
M(\lam)}r_{(S,U,\lam),(Y,U,\lam),(S,U,\lam)}C_{S,U}^{\lam}D_{U,Y}^{\lam}D_{V,S}^{\lam}\\&=&\sum_{Y\in
M(\lam)}\Phi(Y,U)\Psi(Y,V)C_{S,U}^{\lam}D_{U,S}^{\lam}.
\end{eqnarray*}
Then $\sum\limits_{Y\in
M(\lam)}\Phi(Y,U)\Psi(Y,V)C_{S,U}^{\lam}D_{U,S}^{\lam}=0$.
This implies that $$\tau(\sum\limits_{Y\in
M(\lam)}\Phi(Y,U)\Psi(Y,V)C_{S,U}^{\lam}D_{U,S}^{\lam})=0.$$ Since
$\tau(C_{S,U}^{\lam}D_{U,S}^{\lam})=1$, then $\sum\limits_{Y\in M(\lam)}\Phi(Y,U)\Psi(Y,V)=0$.
\end{proof}

\begin{cor} \label{2.20}
Let $A$ be a symmetric cellular algebra over an integral
domain $R$. Then $k_{\lam}=0$ for any $\lam\in\Lambda$ with
$\rad\lam\neq 0$.
\end{cor}

\begin{proof}
Since $|G(\lam)|= 0$ is equivalent to
$\rad\lam\neq 0$, then by Lemma \ref{2.19}, $\rad\lam \neq 0$ implies that $k_{\lam}=0$.
\end{proof}

Using the dual basis, for each $\lam\in\Lambda$, we can also define the cell module $W_{D}(\lam)$.
Then the following lemma is clear.

\begin{lem}\label{2.21}
There is a natural isomorphism of $R$-modules
$$D^{\lam}:W_{D}(\lam)\otimes_{R}i(W_{D}(\lam))\rightarrow R{\rm -span}\{D_{S,T}^{\lam}\mid S,T\in M(\lam)\},$$
defined by $(D_{S},D_{T})\rightarrow D_{S,T}^{\lam}$.
\end{lem}

\bigskip

\section{\bf Radicals of Symmetric Cellular Algebras}
\label{xxsec3}

To study radicals of symmetric cellular algebras, we need the following lemma.

\begin{lem}
\label{lmr2}Let $A$ be a symmetric cellular algebra. Then for any $\lam\in\Lambda$, the elements of the form
$\sum\limits_{S,U\in M(\lam)}r_{SU}C_{S,V}^{\lam}D_{V,U}^{\lam}$ with
$r_{SU}\in R$ make an ideal of $A$.
\end{lem}
\begin{proof}
Denote the set of the elements of the form
$\sum\limits_{S,U\in M(\lam)}r_{SU}C_{S,V}^{\lam}D_{V,U}^{\lam}$ by $I^{\lam}$. Then for any $\eta\in\Lambda$, $P,Q\in M(\eta)$, and
$S,U\in M(\lam)$, we claim that the element
$C_{P,Q}^{\eta}C_{S,V}^{\lam}D_{V,U}^{\lam}\in I^{\lam}$. In
fact, by (C3) of Definition \ref{2.4} and Lemma \ref{2.14} (7),
\begin{eqnarray*}
C_{P,Q}^{\eta}C_{S,V}^{\lam}D_{V,U}^{\lam}&=&\sum_{\epsilon\in\Lambda,
X,Y\in
M(\epsilon)}r_{(P,Q,\eta),(S,V,\lam),(X,Y,\epsilon)}C_{X,Y}^{\epsilon}D_{V,U}^{\lam}\\
&=&\sum_{X\in
M(\lam)}r_{(P,Q,\eta),(S,V,\lam),(X,V\lam)}C_{X,V}^{\lam}D_{V,U}^{\lam}
\end{eqnarray*}
The element $C_{S,V}^{\lam}D_{V,U}^{\lam}C_{P,Q}^{\eta}\in I^{\lam}$ is proved similarly.
\end{proof}

We will denote $\sum\limits_{\lam\in\Lambda, k_\lam=0}I^{\lam}$ by $I^\Lambda$.

Similarly, for each $\lam\in\Lambda$, the elements of the form $\sum\limits_{S,U\in M(\lam)}r_{U,S}D_{U,V}^{\lam}C_{V,S}^{\lam}$ with
$r_{U,S}\in R$ also make an ideal $I_{D}^{\lam}$ of $A$. Denote $\sum\limits_{\lam\in\Lambda, k_\lam=0}I_{D}^{\lam}$ by $I_{D}^\Lambda$.

Define $$I=I^\Lambda+I_{D}^\Lambda$$ and define

$\Lambda_{1}=\{\lam\in\Lambda\mid \rad\lam=0\},$\qquad\qquad\qquad\qquad
$\Lambda_{2}=\Lambda_{0}-\Lambda_{1},$

$\Lambda_{3}=\Lambda-\Lambda_{0},$\qquad\qquad\qquad\qquad\qquad\qquad\,\,\,\,\,\,
$\Lambda_{4}=\{\lam\in\Lambda_{1}\mid k_{\lam}=0\}$.

Now we are in a position to give the main results of this paper.

\begin{thm}\label{thm}
Suppose that $R$ is an integral domain and
that $A$ is a symmetric cellular algebra with a cellular basis
$C=\{C_{S,T}^\lam \mid S,T\in M(\lam),\lam\in\Lambda\}$. Let $\tau$ be a symmetrizing trace on $A$ and let $\{D_{T,S}^\lam \mid S,T\in M(\lam),\lam\in\Lambda\}$ be the dual
basis of $C$ with respect to $\tau$. Then\\
{\rm(1)} $I\subseteq \rad A$, $I^{3}=0$.\\
{\rm(2)} $I$ is independent of the choice of $\tau$.\\
Moreover, if $R$ is a field, then\\
{\rm(3)}
$\dim_{R}I\geq\sum\limits_{\lam\in\Lambda_{2}}(n_{\lam}+\dim_{R}\rad\lam)\dim_{R}L_{\lam}+\sum\limits_{\lam\in\Lambda_{4}}n_{\lam}^{2},$
where $n_{\lam}$
is the number of the elements in $M(\lam)$.\\
{\rm(4)}
$\sum\limits_{\lam\in\Lambda_{2}}(\dim_{K}L_{\lam})^{2}-\sum\limits_{\lam\in\Lambda_{3}}n_{\lam}^{2}\leq
\sum\limits_{\lam\in\Lambda_{2}}(\dim_{K}\rad\lam)^{2}-\sum\limits_{\lam\in\Lambda_{4}}n_{\lam}^{2}.$
\end{thm}

\begin{proof} (1) $I\subseteq \rad A$ , $I^{3}=0$.\\

Firstly, we prove $(I^{\Lambda})^2=0$. Obviously, by the definition of $I^\Lambda$, every element of $(I^{\Lambda})^2$ can be written as a linear combination of elements of the form
$C_{S_{1},T}^{\lam}D_{T,S_{2}}^{\lam}C_{U_{1},V}^{\mu}D_{V,U_{2}}^{\mu}$(we
omit the coefficient here) with $k_\lam=0$ and $k_\mu=0$.

If $\mu<\lam$, then $C_{S_{1},T}^{\lam}D_{T,S_{2}}^{\lam}C_{U_{1},V}^{\mu}D_{V,U_{2}}^{\mu}=0$ by Lemma \ref{2.14} (8).

If $\mu>\lam$, then by Lemma \ref{2.14} (1) and (7),
$$C_{S_{1},T}^{\lam}D_{T,S_{2}}^{\lam}C_{U_{1},V}^{\mu}D_{V,U_{2}}^{\mu}=
\sum_{Y\in
M(\lam)}r_{(U_{1},V,\mu),(Y,T,\lam),(S_{2},T,\lam)}C_{S_{1},T}^{\lam}D_{T,Y}^{\lam}D_{V,U_{2}}^{\mu}.$$
However, by Lemma \ref{2.15}, every $D_{P,Q}^{\eta}$ with nonzero coefficient in the
expansion of $D_{T,Y}^{\lam}D_{V,U_{2}}^{\mu}$ satisfies
$\eta\geq\mu$. Since $\mu>\lam$, then $\eta>\lam$. Now, by Lemma \ref{2.14} (7),
we have $C_{S_{1},T}^{\lam}D_{P,Q}^{\eta}=0$, that is,
$C_{S_{1},T}^{\lam}D_{T,S_{2}}^{\lam}C_{U_{1},V}^{\mu}D_{V,U_{2}}^{\mu}=0$
if $\mu>\lam$.

If $\lam=\mu$, by Lemma \ref{2.14} (3) and (4), we only need to consider the elements of the form
$$C_{S_{1},T_{1}}^{\lam}D_{T_{1},S_{2}}^{\lam}C_{S_{2},T_{2}}^{\lam}D_{T_{2},S_{3}}^{\lam}.$$
By Lemma \ref{2.17} and Lemma \ref{2.20},
\begin{eqnarray*}
C_{S_{1},T_{1}}^{\lam}D_{T_{1},S_{2}}^{\lam}C_{S_{2},T_{2}}^{\lam}D_{T_{2},S_{3}}^{\lam}=
k_{\lam}C_{S_{1},T_{1}}^{\lam}D_{T_{1},S_{3}}^{\lam}=0.
\end{eqnarray*}
Then we get that all the elements of the form
$C_{S_{1},T}^{\lam}D_{T,S_{2}}^{\lam}C_{U_{1},V}^{\mu}D_{V,U_{2}}^{\mu}$
are zero, that is, $(I^{\Lambda})^2=0$.

Similarly, we get $(I_{D}^{\Lambda})^2=0$.

To prove $I^3=0$, we now only need to consider the elements in $I^{\Lambda}I_{D}^{\Lambda}I^{\Lambda}$ and $I_{D}^{\Lambda}I^{\Lambda}I_{D}^{\Lambda}$. For $\lam,\mu,\eta\in\Lambda$ with $k_\lam=k_\mu=k_\eta=0$ and $S,T,M\in M(\lam)$, $U,V,N\in M(\mu)$, $P,Q,W\in M(\eta)$, suppose that $C_{S,T}^{\lam}D_{T,M}^{\lam}D_{U,V}^{\mu}C_{V,N}^{\mu}C_{P,Q}^{\eta}D_{Q,W}^{\eta}\neq 0$. If $\lam>\mu$, then any $D_{X,Y}^{\epsilon}$ with nonzero coefficient in the expansion of $D_{T,M}^{\lam}D_{U,V}^{\mu}$ satisfies $\epsilon\geq\lam$, so $\epsilon>\mu$, this implies that $D_{X,Y}^{\epsilon}C_{V,N}^{\mu}=0$ by Lemma \ref{2.14}, a contradiction. If $\lam<\mu$, then any $D_{X,Y}^{\epsilon}$ with nonzero coefficient in the expansion of $D_{T,M}^{\lam}D_{U,V}^{\mu}$ satisfies $\epsilon\geq\mu$, so $\epsilon>\lam$, this implies that $C_{S,T}^{\lam}D_{X,Y}^{\epsilon}=0$ by Lemma \ref{2.14}, a contradiction. Thus $\lam=\mu$. Similarly, we get $\eta=\mu$. By a direct computation, we can also get $C_{S,T}^{\lam}D_{T,M}^{\lam}D_{U,V}^{\mu}C_{V,N}^{\mu}C_{P,Q}^{\eta}D_{Q,W}^{\eta}=0$. This implies that $I^{\Lambda}I_{D}^{\Lambda}I^{\Lambda}=0$. Similarly $I_{D}^{\Lambda}I^{\Lambda}I_{D}^{\Lambda}=0$ is proved. Then $I^3=0$ follows.

Now it is clear that
$I\subseteq \rad A$ for $I$ is a nilpotent ideal of $A$.\\

(2)\,\, $I$ is independent of the choice of $\tau$.\\

Let $\tau$ and $\tau'$ be two symmetrizing traces and $D$, $d$ the
dual bases determined by $\tau$ and $\tau'$ respectively. By
Lemma \ref{2.3}, for arbitrary $d_{U,V}^{\lam}\in d$,
$$d_{U,V}^{\lam}=\sum_{\varepsilon\in\Lambda, X,Y\in
M(\varepsilon)}\tau(C_{X,Y}^{\varepsilon}d_{U,V}^{\lam})D_{Y,X}^{\varepsilon}.$$
Then for arbitrary $S\in M(\lam)$,
$$C_{S,U}^{\lam}d_{U,V}^{\lam}=\sum_{\varepsilon\in\Lambda, X,Y\in
M(\varepsilon)}\tau(C_{X,Y}^{\varepsilon}d_{U,V}^{\lam})C_{S,U}^{\lam}D_{Y,X}^{\varepsilon}.$$
By Lemma \ref{2.14} (7), (8), if $\varepsilon<\lam$, then $C_{X,Y}^{\varepsilon}d_{U,V}^{\lam}=0$; if $\varepsilon>\lam$, then $C_{S,U}^{\lam}D_{Y,X}^{\varepsilon}=0.$ This implies that
$$C_{S,U}^{\lam}d_{U,V}^{\lam}=\sum_{X,Y\in
M(\lam)}\tau(C_{X,Y}^{\lam}d_{U,V}^{\lam})C_{S,U}^{\lam}D_{Y,X}^{\lam}.$$
By Lemma \ref{2.14} (5), if $Y\neq U$, then $C_{S,U}^{\lam}D_{Y,X}^{\lam}=0$. Hence
$$C_{S,U}^{\lam}d_{U,V}^{\lam}=\sum_{X\in
M(\lam)}\tau(C_{X,U}^{\lam}d_{U,V}^{\lam})C_{S,U}^{\lam}D_{U,X}^{\lam}.$$
Noting that $\tau(C_{X,U}^{\lam}d_{U,V}^{\lam})=\tau(d_{U,V}^{\lam}C_{X,U}^{\lam})$, it follows from Lemma \ref{2.14} that $d_{U,V}^{\lam}C_{X,U}^{\lam}=0$ if $X\neq V$. Thus $$C_{S,U}^{\lam}d_{U,V}^{\lam}=\tau(C_{V,U}^{\lam}d_{U,V}^{\lam})C_{S,U}^{\lam}D_{U,V}^{\lam}.$$

Similarly, we obtain $$C_{S,U}^{\lam}D_{U,V}^{\lam}=\tau'(C_{V,U}^{\lam}D_{U,V}^{\lam})C_{S,U}^{\lam}d_{U,V}^{\lam},$$
$$d_{V,U}^{\lam}C_{U,S}^{\lam}=\tau(C_{V,U}^{\lam}d_{U,V}^{\lam})D_{V,U}^{\lam}C_{U,S}^{\lam},$$
$$D_{V,U}^{\lam}C_{U,S}^{\lam}=\tau'(C_{V,U}^{\lam}D_{U,V}^{\lam})d_{V,U}^{\lam}C_{U,S}^{\lam}.$$
The above four formulas imply that $I$ is independent of the choice
of symmetrizing trace.\\

(3) $\dim_{R}I\geq\sum\limits_{\lam\in\Lambda_{2}}(n_{\lam}+\dim_{R}\rad\lam)\dim_{R}L_{\lam}+\sum\limits_{\lam\in\Lambda_{4}}n_{\lam}^{2}.$\\

For any $\lam\in\Lambda_{2}$ and $S,T\in M(\lam)$, it follows from Lemma \ref{2.14} that
$$C_{S,T}^{\lam}D_{T,T}^{\lam}\equiv \sum\limits_{X\in M(\lam)}\Phi(X,S)D_{X,T}^{\lam}\,\,\,\,\,(\mod A_{D}(>\lam)),$$
$$D_{T,T}^{\lam}C_{T,S}^{\lam}\equiv \sum\limits_{Y\in M(\lam)}\Phi(Y,S)D_{T,Y}^{\lam}\,\,\,\,\,(\mod A_{D}(>\lam)).$$
Let $V$ be the $R$-space generated by
$$\{\sum\limits_{X\in M(\lam)}\Phi(X,S)D_{X,T}^{\lam}\mid S,T\in M(\lam)\}\cup \{\sum\limits_{Y\in M(\lam)}\Phi(Y,S)D_{T,Y}^{\lam}\mid S,T\in M(\lam)\}.$$
Then it is easy to know from the definition of $I^\lam$ and $I_{D}^\lam$ that $$\dim_{R}(I^{\lam}+I_{D}^\lam)\geq\dim V.$$ Note that by Lemma \ref{2.21},
$D^{\lam}: (D_{S},D_{T})\rightarrow D_{S,T}^{\lam}$ is an isomorphism of $R$-modules. So we only need to consider the dimension of $V'$ generated by
$$\{\sum\limits_{X\in M(\lam)}\Phi(X,S)D_{X}\otimes D_{T}\mid S,T\in M(\lam)\}\cup \{D_{T}\otimes\sum\limits_{Y\in M(\lam)}\Phi(Y,S)D_{Y}\mid S,T\in M(\lam)\}.$$ Since $\Phi_{\lam}\neq 0$, $\rank G_{\lam}=\dim_{R}L_{\lam}$, we have $\dim V'=2n_{\lam}\dim_{R}L_{\lam}-(\dim_{R}L_{\lam})^2,$ that is,
$\dim V'=\dim_{R}L_{\lam}\times(n_{\lam}+\dim_{R}\rad\lam)$. Thus $$\dim_{R}(I^{\lam}+I_{D}^{\lam})\geq \dim_{R}L_{\lam}\times(n_{\lam}+\dim_{R}\rad\lam).$$

Clearly, the above inequality holds true for any $\lam\in\Lambda_{4}$, then we have $$\dim_{R}(I^{\lam}+I_{D}^{\lam})\geq n_{\lam}^{2}$$ for any $\lam\in\Lambda_{4}$.

It is clear from Lemma \ref{2.15} that $\dim_{R}I\geq \sum\limits_{\lam\in\Lambda_{2}}\dim_{R}(I^{\lam}+I_{D}^{\lam})+\sum\limits_{\lam\in\Lambda_{4}}n_{\lam}^{2}$ and then item (3) follows.\\

(4)
$\sum\limits_{\lam\in\Lambda_{2}}(\dim_{K}L_{\lam})^{2}-\sum\limits_{\lam\in\Lambda_{3}}n_{\lam}^{2}\leq
\sum\limits_{\lam\in\Lambda_{2}}(\dim_{K}\rad\lam)^{2}.$\\

By (1) and (3),
$$\dim_{R}\rad
A\geq\sum\limits_{\lam\in\Lambda_{2}}(n_{\lam}+\dim_{R}\rad\lam)\dim_{R}L_{\lam}+\sum\limits_{\lam\in\Lambda_{4}}n_{\lam}^{2}.$$ By
the formula
$$\dim_{R}\rad A=\dim_{R}A-\sum_{\lam\in\Lambda_{0}}(\dim_{R}L_{\lam})^{2},$$ we have
$$\dim_{R}A-\sum_{\lam\in\Lambda_{0}}(\dim_{R}L_{\lam})^{2}\geq
\sum\limits_{\lam\in\Lambda_{2}}(n_{\lam}+\dim_{R}\rad\lam)\dim_{R}L_{\lam}+\sum\limits_{\lam\in\Lambda_{4}}n_{\lam}^{2}.$$  That
is,
$$\sum_{\lam\in\Lambda_{3}}n_{\lam}^{2}+\sum_{\lam\in\Lambda_{0}}n_{\lam}^{2}-\sum_{\lam\in\Lambda_{0}}(\dim_{R}L_{\lam})^{2}\geq
\sum\limits_{\lam\in\Lambda_{2}}(n_{\lam}+\dim_{R}\rad\lam)\dim_{R}L_{\lam}+\sum\limits_{\lam\in\Lambda_{4}}n_{\lam}^{2},$$ or
$$\sum_{\lam\in\Lambda_{3}}n_{\lam}^{2}+\sum_{\lam\in\Lambda_{2}}n_{\lam}^{2}-\sum_{\lam\in\Lambda_{2}}(\dim_{R}L_{\lam})^{2}\geq
\sum\limits_{\lam\in\Lambda_{2}}(n_{\lam}+\dim_{R}\rad\lam)\dim_{R}L_{\lam}+\sum\limits_{\lam\in\Lambda_{4}}n_{\lam}^{2},$$
or $$\sum\limits_{\lam\in\Lambda_{2}}(\dim_{K}L_{\lam})^{2}-\sum\limits_{\lam\in\Lambda_{3}}n_{\lam}^{2}\leq \sum\limits_{\lam\in\Lambda_{2}}n_{\lam}^{2}-\sum\limits_{\lam\in\Lambda_{2}}(n_{\lam}+
\dim_{R}\rad\lam)\dim_{R}L_{\lam}-\sum\limits_{\lam\in\Lambda_{4}}n_{\lam}^{2}.$$
According to $\dim_{R}L_{\lam}=n_{\lam}-\dim_{R}\rad\lam$, the right side of the above inequality is $\sum\limits_{\lam\in\Lambda_{2}}(\dim_{K}\rad\lam)^{2}-\sum\limits_{\lam\in\Lambda_{4}}n_{\lam}^{2}$ and this completes the proof.
\end{proof}

\begin{cor} Let $R$ be an integral domain and $A$ a symmetric
cellular algebra. Let $\lam$ be the minimal element in $\Lambda$. If
$\rad\lam\neq 0$, then $R-{\rm span}\{C_{S,T}^{\lam}\mid S,T\in M(\lam)\}\subset \rad A$.
\end{cor}
\begin{proof}
If $a=\sum\limits_{X,Y\in M(\lam)}r_{X,Y}C_{X,Y}^{\lam}$ is not in $\rad A$, then there exists some
$D_{U,V}^{\mu}$ such that $aD_{U,V}^{\mu}\notin \rad A$. If
$\mu\neq\lam$, then $aD_{U,V}^{\mu}=0$ by Lemma \ref{2.14}, it
is in $\rad A$. If $\mu=\lam$, then $aD_{U,V}^{\mu}\in \rad A$ by
Theorem \ref{thm}. It is a contradiction.
\end{proof}

\begin{cor}
Let $A$ be a finite dimensional symmetric cellular algebra and $r\in \rad A$. Assume that $\lam\in\Lambda$ satisfies:\\
{\rm(1)} There exists $S,T\in M(\lam)$ such that $C_{S,T}^{\lam}$ appears in the expansion of $r$ with nonzero coefficient.\\
{\rm(2)} For any $\mu>\lam$ and $U,V\in M(\mu)$, the coefficient of $C_{U,V}^{\mu}$ in the expansion of $r$ is zero.\\
Then $k_{\lam}= 0$.
\end{cor}
\begin{proof} Since $r=\sum\limits_{\varepsilon\in\Lambda, X,Y\in M(\varepsilon)}r_{X,Y,\varepsilon}C_{X,Y}^{\varepsilon}\in \rad A$, we have $rD_{T,S}^{\lam}\in
\rad A$. The conditions (1) and (2) imply that
$$rD_{T,S}^{\lam}=\sum\limits_{X\in
M(\lam)}r_{X,T,\lam}C_{X,T}^{\lam}D_{T,S}^{\lam}.$$It is easy to check that $(rD_{T,S}^{\lam})^{n}=(k_{\lam}r_{S,T,\lam})^{n-1}rD_{T,S}^{\lam}$. Applying $\tau$ on both sides of this equation, we get  $\tau((rD_{T,S}^{\lam})^{n})=(k_{\lam}r_{S,T,\lam})^{n-1}r_{S,T,\lam}$.
If $k_{\lam}\neq
0$, then $\tau((rD_{T,S}^{\lam})^{n})\neq 0$. Hence $rD_{T,S}^{\lam}$ is not nilpotent and then
$rD_{T,S}^{\lam}\notin \rad A$, a contradiction. This implies that
$k_{\lam}= 0$.
\end{proof}

\noindent{\bf Example} The group algebra $\mathbb{Z}_{3}S_{3}$.

The algebra has a basis
$$\{1, s_{1}, s_{2}, s_{1}s_{2}, s_{2}s_{1}, s_{1}s_{2}s_{1}\}.$$
A cellular basis is

$C_{1,1}^{(3)}=1+s_{1}+s_{2}+s_{1}s_{2}+s_{2}s_{1}+s_{1}s_{2}s_{1}$,

$C_{1,1}^{(2,1)}=1+s_{1},
\,\,\,\,\,\,\,\,\,C_{1,2}^{(2,1)}=s_{2}+s_{1}s_{2}$,

$C_{2,1}^{(2,1)}=s_{2}+s_{2}s_{1},
C_{2,2}^{(2,1)}=1+s_{1}s_{2}s_{1}$,

$C_{1,1}^{(1^3)}=1$.\\The corresponding dual basis is

$D_{1,1}^{(3)}=-s_{2}+s_{1}s_{2}+s_{2}s_{1}$,

$D_{1,1}^{(2,1)}=s_{1}+s_{2}-s_{1}s_{2}-s_{2}s_{1},
D_{2,1}^{(2,1)}=s_{2}-s_{1}s_{2}$,

$D_{1,2}^{(2,1)}=s_{2}-s_{2}s_{1},\,\,\,\,\,\,\,\,\,\,\,\,\,\,\,\,\,\,\,\,\,\,\,\,\,\,\,\,\,\,\,\,\,\,\,\,\,\,\,\,\,
D_{2,2}^{(2,1)}=s_{2}-s_{1}s_{2}-s_{2}s_{1}+s_{1}s_{2}s_{1}$,

$D_{1,1}^{(1^3)}=1-s_{1}-s_{2}+s_{1}s_{2}+s_{2}s_{1}-s_{1}s_{2}s_{1}$.

It is easy to know that $\Lambda_{3}=(3)$ and $\Lambda_{1}=(1^3)$.
Then $\dim_{K}\rad A=4$. Now we compute $I$.

$C_{1,1}^{(3)}D_{1,1}^{(3)}=1+s_{1}+s_{2}+s_{1}s_{2}+s_{2}s_{1}+s_{1}s_{2}s_{1}$,

$C_{1,2}^{(2,1)}D_{2,1}^{(2,1)}=1+s_{1}-s_{2}-s_{1}s_{2}s_{1}$,

$C_{1,2}^{(2,1)}D_{2,2}^{(2,1)}=s_{2}+s_{1}s_{2}-s_{2}s_{1}-s_{1}s_{2}s_{1}$,

$C_{2,1}^{(2,1)}D_{1,2}^{(2,1)}=1-s_{1}-s_{1}s_{2}+s_{1}s_{2}s_{1}$,

$C_{2,1}^{(2,1)}D_{1,1}^{(2,1)}=s_{2}+s_{2}s_{1}-s_{1}-s_{1}s_{2}$.

Then $\dim_{K} I=4$. This implies that $I=\rad A$.

\bigskip

\section{\bf Semisimplicity of symmetric cellular algebras}

As a by-product of the results on radicals, we will give some
equivalent conditions for a finite dimensional symmetric cellular
algebra to be semisimple.

\begin{cor}\label{3.5}
Let $A$ be a finite dimensional symmetric cellular
algebra. Then the following are equivalent.\\
{\rm(1)} The algebra $A$ is semisimple.\\
{\rm(2)} $k_{\lam}\neq 0$ for all $\lam\in\Lambda$.\\
{\rm(3)} $\{C_{S,T}^{\lam}D_{T,T}^{\lam}\mid\lam\in\Lambda, S,T\in
M(\lam)\}$ is a basis of $A$.\\
{\rm(4)} For any $\lam\in\Lambda$, there exist $S,T\in M(\lam)$, such
that $(C_{S,T}^{\lam}D_{T,S}^{\lam})^{2}\neq 0$.\\
{\rm(5)} For any $\lam\in\Lambda$ and arbitrary $S,T\in M(\lam)$,
$(C_{S,T}^{\lam}D_{T,S}^{\lam})^{2}\neq 0$.
\end{cor}
\begin{proof}
(2)$\Longrightarrow$(1) If $k_{\lam}\neq 0$ for
all $\lam\in\Lambda$, then $\rad\lam=0$ for all $\lam\in\Lambda$ by
Corollary \ref{2.20}. This implies that $A$ is semisimple by Theorem
\ref{glthm}.

(1)$\Longrightarrow$(2) Assume that there exists some
$\lam\in\Lambda$ such that $k_{\lam}=0$. Then it is easy to check
that $I^{\lam}$ is a nilpotent ideal of $A$. Obviously,
$I^{\lam}\neq 0$ because at least
$C_{U,V}^{\lam}D_{V,U}^{\lam}\neq 0$. This implies that
$I^{\lam}\subseteq \rad A$. But $A$ is semisimple, a
contradiction. This implies that $k_{\lam}\neq 0$ for all
$\lam\in\Lambda$.

(2)$\Longrightarrow$(3) Let $\sum\limits_{\lam\in\Lambda, S,T\in
M(\lam)}k_{S,T,\lam}C_{S,T}^{\lam}D_{T,T}^{\lam}=0$. Take a maximal
element $\lam_{0}\in\Lambda$. For arbitrary $X,Y\in M(\lam_{0})$,
\begin{eqnarray*}
C_{X,X}^{\lam_{0}}D_{X,Y}^{\lam_{0}}(\sum_{\lam\in\Lambda, S,T\in
M(\lam)}k_{S,T,\lam}C_{S,T}^{\lam}D_{T,T}^{\lam})=k_{\lam_{0}}\sum_{T\in
M(\lam_{0})}k_{Y,T,\lam_{0}}C_{X,T}^{\lam_{0}}D_{T,T}^{\lam_{0}}=0.
\end{eqnarray*}
This implies that $\tau(k_{\lam_{0}}\sum\limits_{T\in
M(\lam_{0})}k_{Y,T,\lam_{0}}C_{X,T}^{\lam_{0}}D_{T,T}^{\lam_{0}})=0$,
i.e., $k_{\lam_{0}}k_{Y,X,\lam_{0}}=0$. Since $k_{\lam_{0}}\neq 0$, then we get
$k_{Y,X,\lam_{0}}=0$.

Repeating the process as above, we get that all the $k_{S,T,\lam}$ are
zeros.

(3)$\Longrightarrow$(2) Since
$\{C_{S,T}^{\lam}D_{T,T}^{\lam}\mid\lam\in\Lambda, S,T\in M(\lam)\}$
is a basis of $A$, we have $$1=\sum_{\lam\in\Lambda, S,T\in
M(\lam)}k_{S,T,\lam}C_{S,T}^{\lam}D_{T,T}^{\lam}.$$ For arbitrary
$\mu\in\Lambda$ and $U,V\in M(\mu)$, we have
\begin{eqnarray*}
C_{U,V}^{\mu}D_{V,V}^{\mu}&=&\sum_{\lam\in\Lambda, S,T\in
M(\lam)}k_{S,T,\lam}C_{S,T}^{\lam}D_{T,T}^{\lam}C_{U,V}^{\mu}D_{V,V}^{\mu}\\
&=&k_{\mu}\sum_{ X\in M(\mu)}k_{X,U,\mu}C_{X,V}^{\mu}D_{V,V}^{\mu}.
\end{eqnarray*}
This implies that $k_{\mu}\neq 0$ since
$C_{U,V}^{\mu}D_{V,V}^{\mu}\neq 0$. The fact that $\mu$ is arbitrary
implies that $k_{\lam}\neq 0$ for all $\lam\in\Lambda$.

(2)$\Longleftrightarrow$(4) and (2)$\Longleftrightarrow$(5) are
clear by Lemma \ref{2.17}.
\end{proof}

\begin{cor} Let $R$ be an integral domain and $A$ a symmetric cellular algebra with a cell datum $(\Lambda, M, C, i)$. Let $K$ be the field of fractions of $R$ and $A_{K}=A\bigotimes_{R}K$. If $A_{K}$ is semisimple, then
$$\{\mathcal
{E}_{S,T}^{\lam}=C_{S,S}^{\lam}D_{S,T}^{\lam}C_{T,T}^{\lam}\mid\lam\in\Lambda,S,T\in
M(\lam)\}$$ is a cellular basis of $A_{K}$. Moreover,  if
$\lam\neq\mu$, then $\mathcal {E}_{S,T}^{\lam}\mathcal
{E}_{U,V}^{\mu}=0$.
\end{cor}
\begin{proof} Firstly, we prove that $\{\mathcal
{E}_{S,T}^{\lam}\mid\lam\in\Lambda,S,T\in M(\lam)\}$ is a basis of
$A_{K}$. We only need to show the elements in this set are
$K$-linear independent.
By Lemma \ref{2.14}, we have
\begin{eqnarray*}
\mathcal {E}_{S,T}^{\lam}
&=&\sum\limits_{X\in M(\lam)}r_{(T,T,\lam),(X,S,\lam),(T,S,\lam)}C_{S,S}^{\lam}D_{S,X}^{\lam}\\
&=&\sum\limits_{X\in
M(\lam)}\Phi(X,T)C_{S,X}^{\lam}D_{X,X}^{\lam}
\end{eqnarray*}
for all $\lam\in\Lambda, S,T\in M(\lam)$. Since $A_{K}$ is semisimple,
all $G(\lam)$ are non-degenerate.
Moreover, $\{C_{S,T}^{\lam}D_{T,T}^{\lam}\mid\lam\in\Lambda, S,T\in
M(\lam)\}$ is a basis of $A_{K}$ by Corollary \ref{3.5}, then
$$\{\mathcal
{E}_{S,T}^{\lam}=C_{S,S}^{\lam}D_{S,T}^{\lam}C_{T,T}^{\lam}\mid\lam\in\Lambda,S,T\in
M(\lam)\}$$ is a basis of $A_{K}$.

Secondly, $i(\mathcal {E}_{S,T}^{\lam})\equiv\mathcal {E}_{T,S}^{\lam}$ for arbitrary $\lam\in\Lambda$, and $S,T\in M(\lam)$.
This is clear by Lemma \ref{2.14} and \ref{2.15}.

Thirdly, for arbitrary $a\in A$, since $\{C_{S,T}^{\lam}\mid\lam\in\Lambda, S,T\in M(\lam)\}$ is a cellular basis of $A$, we have
\begin{eqnarray*}
a\mathcal{E}_{S,T}^{\lam}&=&aC_{S,S}^{\lam}D_{S,T}^{\lam}C_{T,T}^{\lam}\\
&=&\sum_{X\in
M(\lam)}r_{a}(X,S)C_{X,S}^{\lam}
D_{S,T}^{\lam}C_{T,T}^{\lam}\\
&=&\sum_{X\in
M(\lam)}r_{a}(X,S)C_{X,X}^{\lam}
D_{X,T}^{\lam}C_{T,T}^{\lam}\\
&=&\sum_{X\in
M(\lam)}r_{a}(X,S)\mathcal
{E}_{X,T}^{\lam}.
\end{eqnarray*}
Clearly, $r_{a}(X,S)$ is
independent of $T$. Then $$\{\mathcal
{E}_{S,T}^{\lam}=C_{S,S}^{\lam}D_{S,T}^{\lam}C_{T,T}^{\lam}\mid\lam\in\Lambda,S,T\in
M(\lam)\}$$ is a cellular basis of $A_{K}$.

Finally, for any $\lam, \mu\in\Lambda$, $S,T\in M(\lam)$, $U,V\in
M(\mu)$,
\begin{eqnarray*}
\mathcal {E}_{S,T}^{\lam}\mathcal
{E}_{U,V}^{\mu}&=&C_{S,S}^{\lam}D_{S,T}^{\lam}C_{T,T}^{\lam}C_{U,U}^{\mu}D_{U,V}^{\mu}C_{V,V}^{\mu}\\
&=&\sum_{\epsilon\in\Lambda, X,Y\in
M(\epsilon)}r_{(T,T,\lam),(U,U,\mu),(X,Y,\epsilon)}
C_{S,S}^{\lam}D_{S,T}^{\lam}C_{X,Y}^{\epsilon}D_{U,V}^{\mu}C_{V,V}^{\mu}.
\end{eqnarray*}
By Lemma \ref{2.14},
$C_{S,S}^{\lam}D_{S,T}^{\lam}C_{X,Y}^{\epsilon}D_{U,V}^{\mu}C_{V,V}^{\mu}\neq
0$ implies $\epsilon\geq\lam, \epsilon\geq\mu$. On the other hand, by Definition
\ref{2.4}, $r_{(T,T,\lam),(U,U,\mu),(X,Y,\epsilon)}\neq 0$ implies
$\epsilon\leq\lam$ and $\epsilon\leq\mu$. Therefore,  if
$\lam\neq\mu$, then $\mathcal {E}_{S,T}^{\lam}\mathcal
{E}_{U,V}^{\mu}=0$.
\end{proof}

\bigskip

\centerline{\bf Acknowledgement}

The author acknowledges his supervisor Prof. C.C. Xi. He also thanks Dr. Wei Hu and Zhankui Xiao for many helpful conversations.


\bigskip


\begin{thebibliography}{}

\bibitem{BS} J. Brundan and C. Stroppel, {\em Highest weight categories arising from Khovanov's diagram algebra I: cellularity}, arxiv: math0806.1532v1.

\bibitem{DR} J. Du and H.B. Rui, {\em Based algebras and standard bases
for quasi-hereditary algebras}, Trans. Amer. Math.Soc., \textbf{350},
(1998), 3207-3235.

\bibitem{G} M. Geck, {\em Hecke algebras of finite type are cellular},
Invent. math., \textbf{169}, (2007), 501-517.

\bibitem{G2} F. Goodman, {\em Cellularity of cyclotomic Birman-Wenzl-Murakami algebras}, J. Algebra, \textbf{321}, (2009), 3299-3320.

\bibitem{G3} J.J. Graham, {\em Modular representations of Hecke algebras and related algebras}, PhD Thesis, Sydney University, 1995.

\bibitem{GL} J.J. Graham and G.I. Lehrer, {\em Cellular algebras},
 Invent. Math., \textbf{123}, (1996), 1-34.

\bibitem{GRM} R.M. Green, {\em Completions of cellular algebras}, Comm.
Algebra, \textbf{27}, (1999), 5349-5366.

\bibitem{GRM2} R.M. Green, {\em Tabular algebras and their asymptotic
versions}, J. Algebra, \textbf{252}, (2002), 27-64.

\bibitem{KL} D. Kazhdan and G. Lusztig, {\em Representations of Coxeter groups
 and Hecke algebras}, Invent. Math., \textbf{53}, (1979), 165-184.

\bibitem{KX1} S. Koenig and C.C. Xi, {\em On the structure of cellular
algebras}, In: I. Reiten, S. Smalo and O. solberg(Eds.): Algebras and
Modules II. Canadian Mathematics Society Proceedings, Vol. \textbf{24},
(1998), 365-386.

\bibitem{KX5} S. Koenig and C.C. Xi, {\em Cellular algebras: Inflations and Morita
equivalences}, J. London Math. Soc. (2), \textbf{60}, (1999), 700-722.

\bibitem{KX7} S. Koenig and C.C. Xi, {\em A characteristic-free approach to Brauer
algebras}, Trans. Amer. Math. Soc., \textbf{353}, (2001), 1489-1505.

\bibitem{KX8} S. Koenig and C.C. Xi, {\em Affine cellular algebras},
preprint.

\bibitem{LZ} G.I. Lehrer and R.B. Zhang, {\em A Temperley-Lieb analogue for the BMW
algebra}, arXiv:math/08060687v1.

\bibitem{MM} G. Malle and A. Mathas, {\em Symmetric cyclotomic Hecke algebras},
J. Algebra, \textbf{205}, (1998), 275-293.

\bibitem{Mu} E. Murphy, {\em The representations of Hecke algebras of type
  $A_{n}$}, J. Algebra, \textbf{173}, (1995), 97-121.

\bibitem{RX} H.B. Rui and C.C. Xi, {\em The representation theory of cyclotomic Temperley-Lieb algebras},
Comment. Math. Helv., \textbf{79}, no.2, (2004), 427-450.

\bibitem{WB} B.W. Westbury, {\em Invariant tensors and cellular categories}, J. Algebra, \textbf{321}, (2009), 3563-3567.

\bibitem{Xi1} C.C. Xi, {\em Partition algebras are cellular}, Compositio
math., \textbf{119}, (1999), 99-109.

\bibitem{Xi2} C.C. Xi, {\em On the quasi-heredity of Birman-Wenzl
algebras}, Adv. Math., \textbf{154}, (2000), 280-298.

\bibitem{XX} C.C. Xi and D.J. Xiang, {\em Cellular algebras and Cartan
matrices}, Linear Algebra Appl., \textbf{365}, (2003), 369-388.

\end{thebibliography}
\end{document}